\documentclass[11pt]{article}

\usepackage{amssymb,amsmath, amsfonts, amsthm}

\usepackage{graphicx}
\usepackage{pgfplots}
\usepackage{tikz}
\usepackage{caption}
\usepackage{booktabs}
\usepackage{array}
\usepackage{verbatim}
\usepackage{subfig}
\usepackage{geometry}
\usepackage{fancyhdr}

\tikzstyle{vertex}=[circle, draw, inner sep=0pt, minimum size=6pt]
\newtheorem{theorem}{Theorem}
\newtheorem{lemma}{Lemma}
\newtheorem{corollary}{Corollary}
\newtheorem{observation}{Observation}

\newtheorem{conjecture}{Conjecture}

\newtheorem{problem}{Problem}

\begin{document}

\title{Well-hued graphs with first difference two}
\author{Geoffrey Boyer$^{a}$ \and Kirsti Kuenzel$^{b}$ \and Jeremy Lyle$^{c}$ \and Ryan Pellico$^{b}$}

\maketitle

\begin{center}
$^a$ School of Mathematical and Statistical Sciences, Clemson University\\
$^b$ Department of Mathematics, Trinity College\\
$^c$ Department of Mathematics and Computer Science, Olivet Nazarene University
\end{center}
\medskip

\maketitle
\begin{abstract}
A graph $G$ is said to be well-hued if every maximal $k$-colorable subgraph of $G$ has the same order $a_k$. Therefore, if $G$ is well-hued, we can associate with $G$ a sequence $\{a_k\}$. Necessary and sufficient conditions were given as to when a sequence $\{a_k\}$ is realized by a well-hued graph. Further, it was conjectured there is only one connected well-hued graph with $a_2 = a_1 + 2$ for every $a_1 \ge 4$. In this paper, we prove this conjecture as well as characterize nearly all well-hued graphs with $a_1=2$. We also investigate when both $G$ and its complement are well-hued. 
\end{abstract}

\maketitle
\section{Introduction}
Throughout this paper, we will consider simple graphs $G$ with vertex set $V(G)$ and edge set $E(G)$. A (proper) $k$-coloring $f$ of a graph $G$ is a mapping $f:\{0, 1, \dots, k\} \rightarrow V(G)$ such that if $(u,v) \in E(G)$, then $f(u) \neq f(v)$.
  A maximal $k$-colorable subgraph is a $k$-colorable subgraph of $G$ induced by a set $S \subseteq V(G)$ such that every graph induced by a set $S'$ with $S\subset S'$ is not $k$-colorable (if $S = V(G)$, then this is a maximal $\chi(G)$-colorable subgraph). 
  
  A graph $G$ is said to be \textit{well-hued} if for any $k \geq 1$, every maximal $k$-colorable subgraph of $G$ has the same order.  We can then associate each well-hued graph $G$ with the sequence $\{a_k\}$ where $a_k$ represents the cardinality of the maximal $k$-colorable subgraph of $G$.  It is apparent that the sequence $a_k$ is non-decreasing and $a_k = n$ for $k \geq \chi(G)$. It is easy to see that any complete graph is a well-hued graph. Another example of a well-hued graph is obtained from the complete graph on six vertices, minus the three edges of a perfect matching.  In this case, $a_1 = 2$, $a_2 = 4$ and $a_k = 6$ for $k\geq 3$.  In particular, every maximal $1$-colorable subgraph of $G$ (i.e. every maximal independent set) has the same order. Therefore, the class of well-hued graphs is a subset of the widely studied well-covered graphs. Plummer \cite{P-1970} first defined well-covered graphs and to date the study of well-covered graphs is concentrated on characterizing large subclasses of them. For instance, well-covered graphs of girth at least $5$ were characterized in \cite{FH-1983, fhn-1993, fhn-1994}. See also the surveys \cite{p-1993} and \cite{H-1999}. In the past couple decades, a wide dirth of similar properties have been studied such as graphs which are well-dominated (all minimal dominating sets have the same order), equimatchable (all maximal matchings have the same order), etc. (See \cite{fhn-1988, LPP-1984}.))
  
   Well-hued graphs were first considered in \cite{GKM-2022}. The authors determined well-hued graphs in many graph families, including cubic and planar graphs, as well as determining conditions for when the join, lexicographic product, or direct product of two graphs is well-hued.  In addition, they characterized the sequences $\{a_k\}$ that could be realized by a well-hued graph, proving the following result.
\begin{theorem}[\cite{GKM-2022}]
  There is a well-hued graph with sequence $\{a_k\}_{k\geq 1}$ if and only if the sequence defined by $d_1 = a_1$ and $d_k =a_k - a_{k-1}$ for $k>1$ is non-increasing and eventually zero.
  \end{theorem}

The graphs that realize the sequences in this theorem are not necessarily connected.  If we restrict our attention to connected graphs, it may no longer be possible to realize a sequence.  In the event that it is possible, there may only be a single graph that realizes the sequence. Indeed, in   \cite{GKM-2022} it was conjectured that there is only one connected well-hued graph with $a_2 = a_1 + 2$ for every $a_1 \geq 4$.  Given a graph $G$, the \textit{corona} of $G$ is the graph obtained by adding a pendant vertex (i.e. leaf vertex) to each vertex of $G$.
\begin{conjecture}[\cite{GKM-2022}] \label{conj:corona}
   For $a_1 \geq 4$, the only connected well-hued graph with $a_2 = a_1 + 2$ is the corona of the complete graph.
  \end{conjecture}

  The goal of this paper is to further examine well-hued graphs.  In Section~\ref{sec:conj} we begin by proving a strengthening of Conjecture \ref{conj:corona};  namely, we prove the following.
  \begin{theorem}
   For $a_1 \geq 3$, the only connected well-hued graph with $a_2 = a_1 + 2$ is the corona of the complete graph.
  \end{theorem}
 If $G$ is a well-hued graph with $a_1 = 1$, it is clear that  $G$ is a complete graph. 
  Therefore, the only well-hued connected graphs $G$  for which $a_2 = a_1 + 2$ that are not the corona of a complete graph are those for which $a_1 = 2$.  In Section~\ref{sec:alpha2}, we nearly characterize all well-hued graphs for which $a_1 = 2$. 
  
  Finally, our investigation of well-hued graphs where $a_1 = 2$ revealed many examples of graphs for which $G$ and its complement, $\overline{G}$, are well-hued.  For this reason, we also consider the class of complement reducible graphs, or cographs; the class of graphs that contains $K_1$ and is closed under the graph operations of disjoint union and join.  In the last section, we determine when cographs are well-hued.

 \subsection{Terminology and Previous Results}
 
 Given a graph $G$, we use $n(G)$ to denote $|V(G)|$. Any vertex in $G$ with degree one is called a \emph{leaf}. A set $I \subseteq V(G)$ is said to be independent if no pair in $I$ are adjacent in $G$. The independence number of $G$, denoted $\alpha(G)$, is the maximum cardinality among all independent sets in $G$. $G$ is said to be \emph{well-covered} if all of its maximal independent sets have the same cardinality, namely $\alpha(G)$. $G$ is bipartite if we can partition $V(G) = A \cup B$ where each of $A$ and $B$ are independent sets in $G$. The disjoint union of graphs $G$ and $H$, denoted $G \cup H$, is the graph such that $V(G \cup H) = V(G) \cup V(H)$ and $E(G \cup H) = E(G) \cup E(H)$; that is, $G\cup H$ is comprised of the disjoint copies of $G$ and $H$.   The join of graphs $G$ and $H$, denoted $G \vee H$, is the graph where 
 $V(G \vee H) = V(G) \cup V(H)$ and $E(G \vee H) = E(G) \cup E(H) \cup (V(G) \times V(H))$; that is, it is obtained from the disjoint union of $G$ and $H$ by adding every  edge between a vertex in $G$ with a vertex in $H$.  The complement of a graph $G$ will be denoted by $\overline{G}$.

 In Section~\ref{sec:alpha2}, we will make use of matchings and alternating paths.  In a graph $G$, a matching $\mathcal{M} \subseteq E(G)$ is a set of independent edges; i.e. edges that do not share any common vertices.  Given a matching $\mathcal{M}$, an alternating path is a subgraph of $G$ that is isomorphic to a path, such that subsequent edges alternate between edges of $\mathcal{M}$ and edges that are not in $\mathcal{M}$.  More specifically, for two matchings $\mathcal{M}$ and $\mathcal{M}'$, we say that a $(\mathcal{M},\mathcal{M}')$ alternating path is a path where edges alternate between edges of $\mathcal{M}$ and $\mathcal{M}'$.  Lastly, given a proper coloring $f$ of a graph $G$, we will call the set of vertices 
$V_k = \{v \in V(G): f(v) = k\}$ the $k$-th color class of the coloring $f$.

The notion of well-hued graphs began as an extension of well-bicovered graphs where $G$ is well-bicovered if every maximal bipartite subgraph of $G$ has the same order (see \cite{GKM-2020}). Thus, if $G$ is well-hued, then it is also well-covered and well-bicovered. 



\section{Graphs with $a_2 = a_1+2$}\label{sec:conj}

\subsection{A general tool for well-covered and well-bicovered graphs}
Here we will introduce a few lemmas to use as tools for graphs that are both well-covered and well-bicovered. The general strategy is to set aside an independent set $I$ in a graph $G$, and for each $v \in V(G) - I$, look at the intersection of the neighborhood of $v$ with $I$.  For that reason, it is helpful to have the following notation. For any maximal independent set $I \subseteq V(G)$,  we define the following:
\begin{itemize}
\item For each $v \in V(G) - I$, let $\Gamma(v,I) = \{u \in V(G)-I: N_G(u) \cap I \subseteq N_G(v) \cap I\}$.
\item Let $I_v$ be any maximum independent set in the graph induced by $\Gamma(v,I)$ containing the vertex $v$.
\end{itemize}

The idea is that we can use the sets $\Gamma(v,I)$ as a relation on the vertices in $V(G) - I$. That is, for any pair of vertices in $V(G)-I$ we define the relation $u\sim v$ if $u \in \Gamma(v, I)$. One can easily verify that this relation is reflexive and transitive. However, we can have vertices $u, v$ such that $N_G(u) \cap I \subset N_G(v) \cap I$ so this relation is not necessarily symmetric. In any case, we will still work with the ``classes" $\Gamma(v,I)$.

\begin{lemma} \label{lem:toolA}
Let $G$ be any graph that is well-covered and let $I \subseteq V(G)$ be any maximal independent set of $G$.  For any vertex $v \in V(G)-I$, $|N_G(v)\cap I| = |I_v|$.
\end{lemma}

\begin{proof}
For graph $G$, maximal independent set $I$, and $v \in V(G)-I$ as indicated in the statement of the lemma, consider the following set of vertices:
\[
I' = (I - (N_G(v) \cap I)) \cup I_v.
\]
This set is a maximal independent set. Indeed, $I'$ is independent as $x \in I_v$ if $N_G(x) \cap I \subseteq N_G(v)\cap I$. Further, $I'$ is maximal in that every vertex in $V(G)-(I \cup \Gamma(v, I))$ is adjacent to a vertex in $I - (N_G(v) \cap I)$ by definition of $\Gamma(v, I)$, every vertex in $\Gamma(v, I) - I_v$ is adjacent to a vertex in $I_v$ (by the maximality of $I_v$), and every vertex in $N_G(v) \cap I$ is adjacent to $v$.
All maximal independent sets have the same order, so $|I| = |I'|$, implying $|N_G(v) \cap I| = |I_v|$.
\end{proof}

This lemma is especially helpful for well-covered and well-bicovered graphs where the difference between the order of a maximal independent set and the maximum order of an induced bipartite subgraph is small. 

\begin{corollary} \label{cor:toolA}
Let $G$ be any graph that is both well-covered and well-bicovered,
 and let $I \subseteq V(G)$ be any maximal independent set of $G$.  If the maximum order of an induced bipartite subgraph of $G$ is $\alpha(G) + k$, then for any $v \in V(G) -I$, $|N_G(v) \cap I| \leq k$.
\end{corollary}
\begin{proof}
The order of a maximal bipartite subgraph is $\alpha(G) + k$.  The set $I \cup I_v$ is the union of two independent sets, and so is clearly bipartite, implying that $k \geq |I_v| = |N_G(v) \cap I|$ (by Lemma \ref{lem:toolA}).
\end{proof}

Furthermore, we can show that if two vertices have neighborhoods that are large enough in $I$, then they must have a common neighbor in $I$.

\begin{lemma} \label{lem:toolB}
Let $G$ be any graph that is both well-covered and well-bicovered,
 and let $I \subseteq V(G)$ be any maximal independent set of $G$.  If the maximum order of an induced bipartite subgraph of $G$ is $\alpha(G) + k$, then for any vertices $x,y \in V(G)-I$, if $|N_G(x)\cap I| + |N_G(y) \cap I| > k$, then $N_G(x) \cap N_G(y) \cap I \neq \emptyset$.
\end{lemma}

\begin{proof}
Choose vertices $x,y \in V(G)-I$. Suppose for contradiction that $N_G(x) \cap N_G(y) \cap I = \emptyset$.  
Then consider the set $I_2 = I \cup I_x \cup I_y$. It is not difficult to see that with the partition 
$A = (N_G(y) \cap I) \cup I_x$ and 
$B = I_y \cup (I - (N_G(y) \cap I))$, 
the graph induced by $I_2$ is bipartite.  Furthermore, since the neighborhoods of $x$ and $y$ in $I$ are disjoint, then $\Gamma(x,I)$ and $\Gamma(y,I)$ are disjoint as well.  
Thus, by Lemma~\ref{lem:toolA},  \[|I_2| = |I| + |I_x| + |I_y| = |I| + |N_G(x)\cap I| + |N_G(y)\cap I| > |I| + k, \] which is a contradiction. Therefore, $N_G(x) \cap N_G(y) \cap I \neq \emptyset$.
\end{proof}

One last lemma allows us to consider the adjacencies inside the graph induced by $V(G)-I$.

\begin{lemma} \label{lem:toolC}
Let $G$ be any graph that is both well-covered and well-bicovered,
 and $I \subseteq V(G)$ be any maximal independent set of $G$.  If the maximum order of an induced bipartite subgraph of $G$ is $\alpha(G) + k$, then for any pair of vertices $x,y \in V(G) - I$, if $|(N_G(x)\cup N_G(y)) \cap I| > k$, then $xy\in E(G)$.
\end{lemma}

\begin{proof}
Choose vertices $x,y \in V(G) - I$ such that $|(N_G(x)\cup N_G(y)) \cap I| = \ell > k$.  (Note that by Corollary \ref{cor:toolA}, we know that $x \not \in \Gamma(y,I)$ and $y \not \in \Gamma(x,I)$). Suppose for the sake of contradiction that $xy\not \in E(G)$ and consider the set 
\[I' = (I - ((N_G(x)\cup N_G(y)) \cap I) \cup \{x,y\}. \]
This set is certainly independent, but the order of the set is $\alpha(G) - \ell + 2$. Since $G$ is well-covered, there exists a set $Z \subset V(G) - I$ of at least $\ell - 2$ vertices such  that  $Z \cup \{x,y\}$ is an independent set. However, this implies  $I \cup Z \cup \{x, y\}$ induces  a bipartite subgraph of order $\alpha(G) + \ell$,  which is a contradiction.  Therefore, $xy \in E(G)$.
\end{proof}

As a quick example of the use of these lemmas, we can provide an alternative proof of a result that originally appeared in \cite{GKM-2022}.
\begin{corollary}[\cite{GKM-2022}] Let $G$ be a well-hued graph without isolated vertices and with $a_2 = a_1 + 1$.  Then $G$ is a complete graph.
\end{corollary}

\begin{proof}
Let $G$ be a well-hued graph with $a_2 = a_1 + 1$, and $I\subseteq V(G)$ be any maximal independent set.  By Corollary \ref{cor:toolA} and Lemma \ref{lem:toolA}, each vertex $v \in V(G) - I$ has exactly one neighbor in $I$, and $\Gamma(v, I)$ forms a clique. By Lemma \ref{lem:toolB}, the neighbor (in $I$) of any pair $u,v \in V(G) - I$ must in fact be the same vertex, so $G$ forms a clique.
\end{proof}

A further simple consequence of these tools.

\begin{corollary}
No well-hued graph contains a vertex adjacent to two leaves.
\end{corollary}

\begin{proof}
Let $G$ be a well-hued graph, and suppose for contradiction that there is a vertex $w\in V(G)$ adjacent to two leaves, say $a$ and $b$.
Then let $I$ be any maximal independent set containing $w$.  By Lemma \ref{lem:toolA}, $|I_a| = |I_b| = 1$, which would imply $a$ is adjacent to $b$, a contradiction.
\end{proof}

Taken together, these tools help to give insight into the structure of well-hued graphs when $a_1$ is large relative to $j = a_2 - a_1$.  In this case, each well-hued graph consists of a union of an independent set $I$, and a dense graph induced by $V(G)-I$ (with maximum independent set of size $a_2- a_1$). Note that we know if $G$ is well-hued, then $a_2 \le 2a_1$ and we know there are connected well-hued graphs with $a_2 = a_1 +j$ for $j \in \{1, 2, a_1\}$. To see that we can find a connected well-hued graph with $a_2 = 2a_1 -1$, consider the graph $H$ obtained from the corona of $C_n$ for odd $n\ge 3$ and attach a $K_2$ to one of the vertices on $C_n$. Thus, $H$ is $3$-colorable, all maximal independent sets have order $n+1$, and all maximal bipartite subgraphs have order $n(H)-1 = 2n+1$. Thus, $a_2 = a_1 + n = 2a_1 -1$. On the other hand, if we replace $C_n$ with $K_n$ in the above example, then $a_1 = n+1$ and $a_2 = n+4 = a_1 + 3$. We would like to know whether there exists a connected well-hued graph with $a_2 = a_1 + j$ for any $4 \le j \le a_1 - 2$?

\subsection{Sequences with first difference equal to two and $\alpha(G)\ge 3$}

Using the results from the previous section, we can show that in fact, the corona of a complete graph is the only well-hued graph with $a_2 = a_1 + 2$ for $a_1 \geq 3$.  This was first conjectured in \cite{GKM-2022} (though with $a_1 \geq 4$), and we show this in the following sequence of results.

\begin{lemma} \label{claim:Parts}
Let $G$ be a well-hued graph with $a_2 = a_1 + 2$ and $a_1 \geq 3$, and let $I$ be any maximal independent set of $G$. 
If there is a vertex $w \in V(G)-I$ such that $|N_G(w) \cap I| = 2$, then the following statements are true:
\begin{enumerate}
\item Every vertex $x \in (V(G) -I )-   \Gamma(w,I)$ satisfies $|N_G(x) \cap I| = 2$.  
\item For any vertex $z \in V(G)-I$ such that $|N_G(z) \cap I| = 1$, then $z \in \Gamma(x,I)$ for every vertex $x \in V(G) - I$.
\item Any pair of vertices $x,y \in V(G)-I$ satisfies $|N_G(x) \cap N_G(y) \cap I| \neq \emptyset$.
\end{enumerate}
\end{lemma}

\begin{proof}
Let $w \in V(G) - I$ be a vertex such that $|N_G(w) \cap I| = 2$.   Then consider another vertex $x \in (V(G)-I) - \Gamma(w, I)$ (such a vertex must exist, as $|I|=a_1 \geq 3$ and $G$ is connected.) Note that $|(N_G(w) \cup N_G(x))\cap I| \ge 3>2$ so by Lemma~\ref{lem:toolC}, $wx \in E(G)$. Further, by Lemma \ref{lem:toolB}, $N_G(w) \cap N_G(x) \cap I \neq \emptyset$, establishing part (1) of the claim.

 For part (2), we note that part (1) implies the only vertices $z \in V(G)-I$ such that $|N_G(z) \cap I| = 1$ must be in $\Gamma(w, I)$.  However, using the argument above, substituting in any vertex in $x \in (V(G)-I) - \Gamma(w,I)$ for $w$, we get that $z$ must also be in $\Gamma(x,I)$. Furthermore, $N_G(x) \cap I$ and $N_G(w) \cap I$ intersect in exactly one vertex  (and again, since $a_1 \geq 3$, there must be a vertex $x \in  (V(G)-I) - \Gamma(w, I)$)    ; this vertex must be 
$N_G(z) \cap I$, so if $y$ is another vertex such that  $|N_G(y) \cap I| = 1$, it must be that $z \in \Gamma(y,I)$ as well (and this establishes part (2) of the claim).

To establish part (3), there are two simple cases to consider for $x, y \in V(G)-I$. The first case is when $|N_G(x) \cap I| = |N_G(y) \cap I| = 2$. In this case, we can simply use the argument above, replacing $w$ with either $x$ or $y$. The second case is when $|N_G(x) \cap I| = 1$ or $|N_G(y) \cap I| = 1$. In this case, by part (2), either $x \in \Gamma(y,I)$ or $y \in \Gamma(x,I)$ implying that $|N_G(x) \cap N_G(y) \cap I| \neq \emptyset$.
\end{proof}

Obviously, if the statements of the above result are not true, then it must be the case that there was no vertex $w \in V(G)-I$ with two neighbors in $I$, which gives the following corollary.

\begin{corollary} \label{cor:Disjoint}
Let $G$ be a well-hued graph with $a_2 = a_1 + 2$ and $a_1 \geq 3$, with maximal independent set $I$. If there are two vertices $x, y \in V(G)-I$ such that $N_G(x) \cap N_G(y) \cap I = \emptyset$, then every vertex in $V(G)-I$ has only one neighbor in $I$.
\end{corollary}

The next step is to show that if there is a vertex $w \in V(G)-I$ such that $|N_G(w) \cap I| = 2$, then we can choose a new independent set $I'$ such that every vertex in $V(G)-I'$ now only has one neighbor in $I'$.


\begin{lemma} \label{claim:DegreeOne}
Let $G$ be a well-hued graph with $a_2 = a_1 + 2$ with $a_1 \geq 3$.
There is a choice of independent set $I$ such that for all vertices $v \in V(G)-I$,  $|N_G(v) \cap I| = 1$.
\end{lemma}

\begin{proof}
If every vertex $w$ satisfies $|N_G(w) \cap I| = 1$, the result follows.  Therefore, we assume that there is a vertex $w$ such that $|N_G(w) \cap I| = 2$. We show that every vertex in $V(G)-I$ is adjacent to the same vertex in $I$. 
This is easy if there is a vertex $z \in V(G)-I$ such that $|N_G(z) \cap I| = 1$.  In this case, set $N_G(z) \cap I = \{w'\}$.  Then, part (3) of Lemma \ref{claim:Parts} implies that every vertex $x \in V(G)-I$ is also adjacent to $w'$.  

Therefore, we may assume that there is no vertex $z \in V(G)-I$ such that $|N_G(z) \cap I| = 1$.  In this case, consider two vertices $x, y \in V(G)-I$ so that $N_G(x) \cap I \neq N_G(y) \cap I$. By part (3) of Lemma \ref{claim:Parts}, we can let $\{a\} = N_G(x) \cap N_G(y) \cap I$ and $N_G(x) \cap I = \{a, b\}$. Then consider the sets $I_x$ and $I_y$.  Since $|N_G(x) \cap I| = |N_G(y) \cap I| = 2$, then $|I_x| = |I_y| = 2$ by Lemma \ref{lem:toolA}.  Furthermore, since there is no vertex $z \in V(G)-I$ such that $N_G(z) \cap I = \{a\}$, the sets $\Gamma(x,I)$ and $\Gamma(y,I)$ are disjoint (and so, in particular, $I_x$ and $I_y$ are disjoint).

Letting $A = I_x \cup (I - N_G(x))$ and $B = I_y \cup \{b\}$, we can see that the set $I_2 = I_x \cup I_y \cup (I- \{a\})$ is a maximal bipartite set of order $a_1 - 1 + 4$, but $a_2 = a_1 + 2$, a contradiction.
Therefore, there is a vertex $w' \in I$ such that every vertex in $V(G)-I$ is adjacent to $w'$. 

Now consider a vertex $w \in V(G)-I$ such that $|N_G(w) \cap I| = 2$,  and a vertex $x \in (V(G)-I) - \Gamma(w,I)$. We consider the set $I_3 = (I -\{w'\})\cup I_w \cup I_x$. We can form the partition 
$A = I_w \cup ((N_G(x) \cap I) - \{w'\} ))$ and 
$B = I_x \cup (I - N_G(x))$, and see that this induces a bipartite graph.  Furthermore, $|I_3|  \geq a_1 + |I_w \cup I_x| - 1$.  Since $|I_w| = |I_x| = 2$, this implies that $I_w$ and $I_x$ have a common intersection, say $w''$.  

Since $w'' \in \Gamma(w,I)$ and $\Gamma(x,I)$, then it must be the case that $|N_G(w'') \cap I| = 1$.  This same argument can be applied to all vertices $x$ such that $|N_G(x) \cap I| = 2$, implying that $w''$ is not adjacent to any vertex in $(V(G)-I)- \Gamma(w,I)$.

Now form a new independent set:
\[
I' = (I - \{w'\}) \cup \{w''\}.
\]
Then $x, w' \in (V(G)-I) - \Gamma(w,I)$ satisfy $|N_G(x) \cap I'| = 1$ and $N_G(w') \cap I' = \{w''\} \ne N_G(x) \cap I'$.  Using Corollary \ref{cor:Disjoint}, this would imply that there is no vertex in $V(G)-I'$ such that $|N_G(w) \cap I'| = 2$.

\end{proof}

Now we are ready to prove the main result. 

\begin{theorem}\label{thm:32}
For any well-hued graph $G$ with $a_2 = a_1 + 2$ and $a_1 \geq 3$, $G$ is the corona of a complete graph.
\end{theorem}

\begin{proof}

Let $G$ be a well-hued graph such that $a_2 = a_1 + 2$ and $a_1 \geq 3$. Let $I \subseteq V(G)$ be a maximal independent set. 
At this point, we simply need to show that $H(I)$, the graph induced by $V(G)-I$,  is a complete graph, and every vertex in $I$ is adjacent to exactly one vertex in $H(I)$ to show that the resulting graph is the corona of a complete graph.

By Lemma~\ref{claim:DegreeOne}, we can choose $I$ such that all vertices in $H(I)$ are adjacent to only one vertex in $I$. Suppose that there are two vertices $x$ and $y$ in $H(I)$ such that $x$ and $y$ are not adjacent. Note by Lemma~\ref{lem:toolA} since $|N_G(x) \cap I| = 1 = |N_G(y)\cap I|$, then $|I_x| = |I_y|=1$. Let $N_G(x) \cap I = \{a\}$ and $N_G(y)\cap I = \{b\}$. This implies that $x$ is adjacent to every vertex in $\Gamma(x, I)$ and $y$ is adjacent to every vertex in $\Gamma(y, I)$. Thus, $x \not\in \Gamma(y, I)$ and we may conclude that $a \ne b$.   Let $ z\in V(G)-I$ be any vertex adjacent to some $c \in I - \{a, b\}$.
Let $I_2 = I\cup  \{x, y, z \}$.
We can partition this as 
$A = \{a,b,z\} \cup (I - \{a,b,c\})$ and 
$B = \{x,y,c\}$.
This would form a bipartite subgraph that is too large.

Now we can finally show that $G$ is the corona of the graph $H(I)$ (for a suitably chosen independent set $I$), which we just showed was a complete graph. Each vertex in $H(I)$ is adjacent to exactly one vertex in $I$.  The last step is to argue that each vertex in $I$ is adjacent to exactly one vertex in $H(I)$.  If this is not the case, let $a$ be a vertex in $I$ with neighbors $x$ and $y$ in $H(I)$.  We can form a bipartite subgraph by choosing $I - \{a\}$ along with $x$ and $y$.  By the claim above, no other vertex from $H(I)$ could be chosen (or it would form a triangle with $x$ and $y$.  Therefore, this is a maximal bipartite subgraph with exactly $a_1 + 1$ vertices, which is too small.  Therefore $G$ is the corona of a complete graph.

\end{proof}

Finally, note that if we associate a difference sequence where each term is of the form $a_k-a_{k-1}$ for the well-hued graph $G$, then the corona of $K_{k-1}$ has the difference sequence $2^11^k$, i.e. the sequence comprised of one $2$, follows by $k$ 1s. Thus,  the result above indicates that sequences like $(4,6,7,8, \dots)$ are realizable (in that there is a well-hued graph with $a_k = 2k+a_1$), but $(4,6,8,9, \dots)$ is not. Indeed,  if there exists some well-hued graph $G$  with sequence $(4,6,8,9, \dots)$, then $a_1 =4$ and $a_2 = a_1+2$, implying $G$ is a corona of a complete graph, contradicting the assumption that $a_3 = a_2+2$. When $a_1 = 2$, there is a nice (and simple) family that demonstrates that every possible difference sequence (that starts with 2) is realizable. In fact, in the next section, we show that the difference sequence $2^{k}1^{\ell}$ is realizable for a graph $G$ with $a_1 = 2$.

\section{Graphs with $\alpha(G)=2$}\label{sec:alpha2}
In this section, we characterize nearly all connected, well-hued graphs with $\alpha(G)=2$. Note that in light of Theorem~\ref{thm:32}, this will nearly characterize all well-hued graphs with first difference $2$. We first attempted a computer search on the NetworkX graph atlas which contains exactly 1252 connected graphs of order $7$ or less. We believe that the following list is a (near) complete list of all $77$ nontrivial, connected well-hued graphs up to order $7$, listed by order, then size, then (lexicographically sorted) degree sequence.

The first part of the title is their position in the atlas, the second is the sequence $[a_1,a_2,\dots]$. The last entry is always $a_{\chi(G)} =n(G)$.
\begin{figure}[h!]
\begin{center}
   \includegraphics[width = 0.85\textwidth]{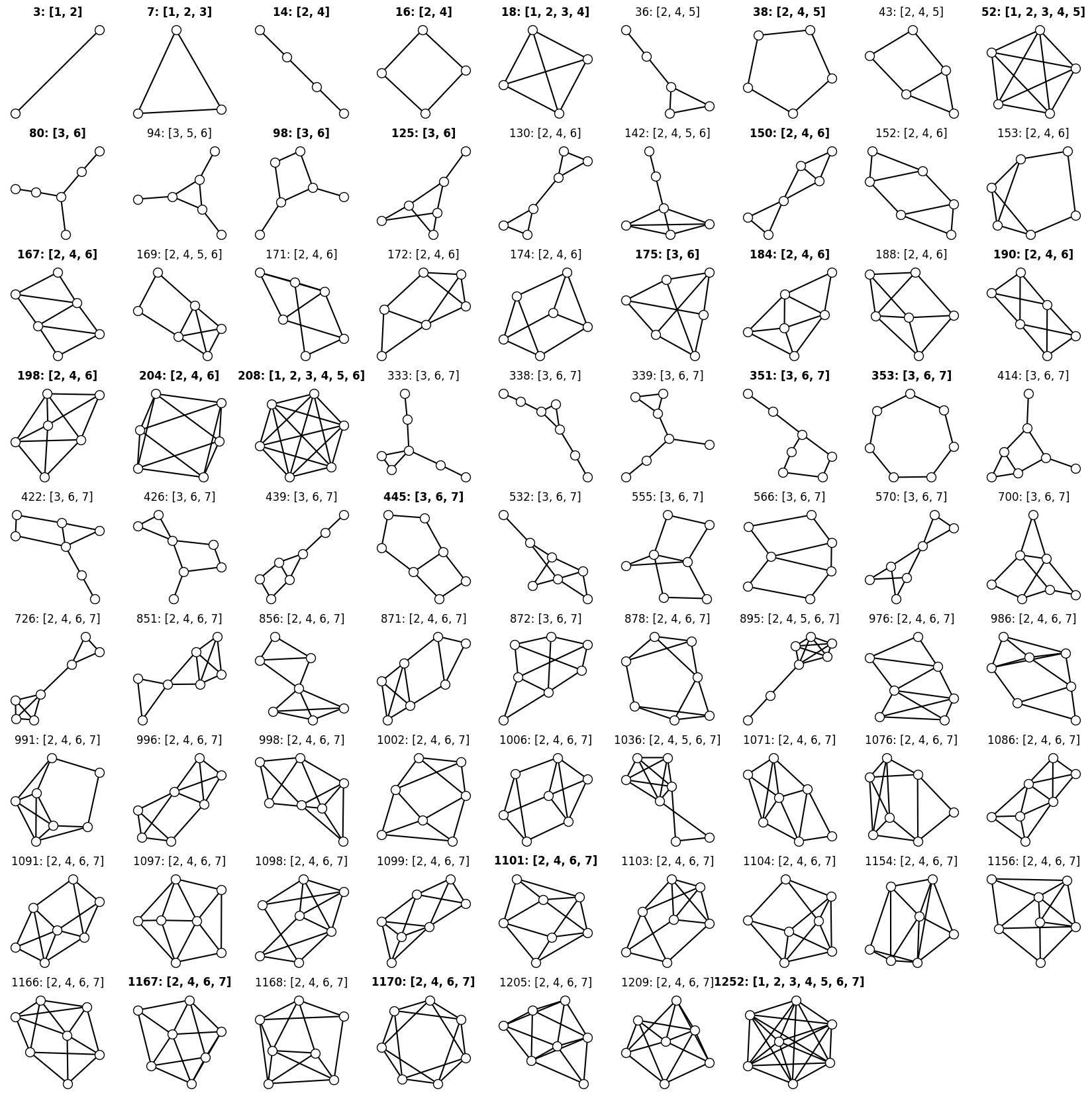}
\end{center}
\caption{Table of well-hued graphs of order at most $7$}
\label{fig:table}
\end{figure}
Among these graphs $G$ we note that all but four (atlas numbers $38, 351, 353, 445$) have the following property: $V(G)$ can be partitioned so that each part $X$ has $|X|\geq 2$ and $G[X]$ is a clique. This observation also aligns with Theorem 2 in the previous section. We also observed the fact that there are many instances when both $G$ and $\overline{G}$ are well-hued. In fact, in Figure~\ref{fig:table}, the following atlas numbers (bolded) represent well-hued graphs whose complement is also well-hued: 3, 7, 14, 16, 18, 38, 52, 80, 98, 125, 150, 167, 175, 184, 190, 198, 204, 208, 351, 353, 445, 1101, 1167, 1170, and 1252. In Section~\ref{sec:cographs}, we will focus on well-hued cographs $G$ and we will see that $G$ is well-hued if and only if  $\overline{G}$ is well-hued. In this section, we focus on well-hued graphs with $\alpha(G)=2$.

First, we characterize the well-hued graphs with $\alpha(G)=2$ and sequence $(2, 4, 6, \dots, 2k)$ where $2k = n(G)$.

\begin{theorem}\label{thm:222} $G$ is a connected, well-hued graph with sequence $(2, 4, 6, \dots, 2k)$ and chromatic number $k$  if and only if it is a spanning subgraph of $K_{2, \dots, 2}$ such that for any $u \in V(G)$, the graph induced by $V(G) - N_G[u]$ is $\overline{K_2}$-free. 
\end{theorem}

\begin{proof}
One direction is clear.  In the other direction, suppose $G$ is a spanning subgraph of $K_{2, \dots, 2}$ such that for any $u \in V(G)$, the graph induced by $V(G) - N_G[u]$ is $\overline{K_2}$-free. We show that $G$ is well-hued. Certainly it is well-covered by the second condition. Let $V_1, \dots, V_k$ be the partite sets of $G$. Let $2 \le j \le k-1$ and consider all $j$-colorable subgraphs of $G$. We know that there are $j$-colorable subgraphs of order $2j$ (namely, take any $j$ partite sets). Note that there are no $j$-colorable subgraphs of order $2j+1$ for this would imply that $G$ contains an independent set of size $3$. So we only need to show that there is not a maximal $j$-colorable subgraph with order at most $2j-1$. To that end, suppose $H$ is such a subgraph and let $W_1, \dots, W_j$ be the color classes of $H$ such that $|W_i|=2$ for $i \in [r]$ and $|W_i|=1$ for $r+1 \le i \le j$.  Let $F = V(G) - V(H)$. We may assume that $W_j$ dominates $F$ for otherwise $H$ is not maximal. However, the vertex in $W_j$ is contained in some $V_{\alpha_j} = \{x_{\alpha_j}, y_{\alpha_j}\}$. We may assume that $W_j = \{x_{\alpha_j}\}$. Thus, $y_{\alpha_j} \in V(H)$. If $y_{\alpha_j}$ is the only vertex in its color class in $H$, then $H$ does not have chromatic number $j$. Thus, reindexing if necessary, we may assume that $y_{\alpha_j} \in W_r$. Let $w$ be the other vertex in $W_r$. If $w$ does not dominate $F$, then there exists some $z \in F$ such that $\{w,z\}$ is an independent set and we could extend $H$ to include $z$ by recoloring $y_{\alpha_j}$ to have color $j$ and color $z$ with $r$. Therefore, we may assume that $w$ dominates $F$. Continuing this same argument, we may assume that $w \in V_{\alpha_{j-1}} = \{x_{\alpha_{j-1}}, y_{\alpha_{j-1}}\}$. With no loss of generality,  assume $w = x_{\alpha_{j-1}}$. It follows that $y_{\alpha_{j-1}} \in V(H)$ and so on until we inevitably reach the contradiction that $H$ is not maximal.  It follows that no such subgraph $H$ exists. Thus, all $j$-colorable subgraphs have order $2j$. 
\end{proof}

Before we move on, we pose the following conjecture about well-hued graphs with sequence $(\alpha(G), 2\alpha(G), 3\alpha(G), \dots, k\alpha(G))$.

\begin{conjecture} $G$ is a connected, well-hued graph with sequence \[(\alpha(G), 2\alpha(G), 3\alpha(G), \dots, k\alpha(G))\] and chromatic number $k$  if and only if it is a spanning subgraph of $K_{\alpha(G), \dots, \alpha(G)}$ such that for any $u \in V(G)$, the graph induced by $V(G) - N_G[u]$ is $\overline{K_{\alpha(G)}}$-free. 
\end{conjecture}

Next, we characterize those well-hued graphs where $\alpha(G)=2$ with chromatic number $k$ and order $2k+1$. To that end, we utilize the idea of alternating paths. Recall that given a matching $\mathcal{M}$ in $G$, we call a $uv$-path $P$ an $(\mathcal{M}, E(G)-\mathcal{M})$-alternating path if no adjacent edges of $P$ are in $\mathcal{M}$ and no adjacent edges of  $P$ are in $E(G)-\mathcal{M}$. Consider a graph $G$ with perfect matching $\mathcal{M}$ (such that $\overline{G}$ is a spanning subgraph of $K_{2,2, \dots, 2}$).  We call a set  $S \subset V(G)$ an \textit{ $\mathcal{M}$-alternating dominating set}  if  every vertex $v \in V(G)$ is either in the set $S$, or lies on an  $(\mathcal{M}, E(G) - \mathcal{M})$-alternating path originating at a vertex $u \in S$ and beginning with the edge from $\mathcal{M}$ involving $u$, i.e. every vertex of $G$ is reachable from $S$ via alternating paths.


First, we note that the matching chosen is irrelevant.

\begin{lemma}
Suppose $G$ has perfect matchings $\mathcal{M}$ and $\mathcal{N}$.  Then $S \subset V(G)$ is an $\mathcal{M}$-alternating dominating set if and only if it is an $\mathcal{N}$-alternating dominating set.
\end{lemma}

\begin{proof}
Let $S$ be an $\mathcal{M}$-alternating dominating set. It is a well-known result that the symmetric difference between perfect matchings $\mathcal{M}$ and $\mathcal{N}$ will consist of $(\mathcal{M}, \mathcal{N})$-alternating even cycles. Therefore, any $(\mathcal{M},E(G) - \mathcal{M})$-alternating path from $u$ to $v$ can be modified to form an $(\mathcal{N},E(G) - \mathcal{N})$-alternating path from $u$ to $v$, and vice versa, simply by traversing any part of the path that intersects the symmetric difference by traveling the opposite direction around these even cycles.

%
\end{proof}

In particular, we will be interested in finding $\mathcal{M}$-alternating dominating sets of a graph $G$ that are also independent, so it is helpful to know that every graph $G$ with a perfect matching has such a set.

\begin{lemma}
Every graph $G$ with perfect matching $\mathcal{M}$ contains an $\mathcal{M}$-alternating dominating set that  is independent.
\end{lemma}

\begin{proof}
Choose $S_0 = \{v_0\} \subset V(G)$ such that $v_0$ reaches the largest number of vertices via an $(\mathcal{M},E(G) - \mathcal{M})$-alternating path (beginning at $v_0$ and starting with the matched edge at $v_0$).  Then, at each additional iteration, if $S_{i-1} = \{v_0, v_1, \dots, v_{i-1}\}$ is not a dominating set, choose a vertex $v_i$ that is not reachable from $S_{i-1}$ which maximizes the number of vertices reachable by $v_i$ that are not reachable by $S$, and add the vertex $v_i$ to $S_i = S_{i-1} \cup \{v_i\}$.

Choosing the set $S$ in this way will ensure that the vertices of $S$ form an independent set in the original graph $G$. Indeed, let the vertices matched up to $S =  \{v_0, v_1, \dots, v_k\}$ in $\mathcal{M}$ be  $\{u_0, u_1, \dots, u_k\}$ respectively.  Suppose that $S$ is not independent in $G$, that is, there is an edge $v_iv_j$ for $i<j$. Note that $u_j$ is not reachable by alternating path from  $S_{i} = \{v_0, v_1, \dots, v_i\}$ in $G$, or else $v_j$ would also be reachable by one of these vertices in $G$.  Furthermore, beginning at $u_j$, any vertex reachable by an alternating path from $v_i$ would also be reachable by an alternating path from $u_j$ (and $v_j$ would also be reachable), which contradicts the choice of $v_i$ as the vertex not reachable by $S_{i-1}$ that maximizes the number of reachable vertices that are not reachable by $S_{i-1}$. 
\end{proof}

In Theorem~\ref{thm:222}, we used the idea of successively recoloring vertices. In that argument, attempting to recolor a vertex has a domino effect (moving one vertex into a color class may force another vertex out of the color class) until we can resolve the coloring (the result is that a $j$-colorable subgraph of order less than $2j$ was not maximal). A coloring with two vertices in each color class can be viewed as a matching in the complement of $G$; alternating paths give an alternative way to view the process of successive recoloring.

\begin{theorem}\label{thm:2k+1}
A connected graph $G$ is well-hued with sequence $(2, 4, 6, \dots, 2k, 2k+1)$ if and only if $G$ is formed in the following way:
\begin{enumerate}
\item 
The graph $G$ is obtained from a vertex $c$ together with a spanning subgraph of $K_{2,\dots, 2}$, call it $G_1$, with order $2k$ and $\alpha(G_1) = 2$. 

\item 
The set $S =  V(G_1) - N_G[c] \neq \emptyset$ is an independent $\mathcal{M}$-alternating dominating set in $\overline{G_1}$ where $\mathcal{M}$ is the matching in $\overline{G_1}$ taken from the partite sets in $G_1$. 
\end{enumerate}

\end{theorem}

\begin{proof}
First, suppose that $G$ is well-hued with sequence  $(2,4, 6, \dots, 2k, 2k+1)$. It is clear that $\alpha(G_1)=2$ and $\chi(G) = n - k = k+1$. 
Consider any coloring of $G$ with $k + 1$ colors.  Since $\alpha(G) = 2$ all color classes must contain either one or two vertices. In fact there must be exactly $k$ color classes with two vertices, so let $V_i = \{x_i, y_i\}$ denote the vertices of color $i\in [k]$.  Furthermore, let $c$ be the vertex in a color class by itself. We define $G_1$ to be the graph induced by the sets $V_i$, i.e. $G_1 = G - \{c\}$. Let $S = V(G_1) - N_G[c]$ and note that $S$ induces  a clique in $G$ for otherwise $\alpha(G)>2$. Thus, $S$ is an independent set in $\overline{G_1}$. All that remains is to show that $S$ is a $\mathcal{M}$-alternating dominating set in $\overline{G_1}$. Suppose that for some $i\in [k]$ that both vertices in $V_i$ are adjacent to $c$ in $G$. Consider $T_i = V(G) - V_i$, which has order $2k-1$ and induces a $k$-colorable subgraph of $G$.  Since $G$ is well-hued with $\alpha_k(G) = 2k$, $T_i$ must not be maximal in $G$, i.e. there is a vertex $w_i \in V_i$ such that $T_i \cup  \{w_i\}$ is also $k$-colorable.

Let $\mathcal{U} = \{U_1, U_2, \dots, U_{k}\}$ be the color classes of $T_i \cup\{w_i\}$,  and let $\mathcal{M}_\mathcal{U}$ denote the corresponding matching in $\overline{G}$.  Then consider the symmetric difference of $\mathcal{M}_\mathcal{U}$ and $\mathcal{M}_\mathcal{V}$ in $\overline{G}$, which consists of $(\mathcal{M}_\mathcal{U}, \mathcal{M}_\mathcal{V})$-alternating paths and even cycles.  Since every vertex from $G_1$ is used in $\mathcal{M}_\mathcal{V}$, then there is exactly one alternating path that begins at $c$ (in $\mathcal{M}_\mathcal{U}$) and ends at the vertex in $V_i - \{w_i\}$.  In this way, we see that both vertices of $V_i$ are reachable along an $(\mathcal{M}_\mathcal{V}, E(\overline{G_1}) - \mathcal{M}_\mathcal{V})$-alternating path of $\overline{G_1}$ from a vertex in $S$.  Repeating this process, we get that $S$ is an $\mathcal{M}$-alternating dominating set of $\overline{G_1}$.\\

Now suppose that $G$ is formed by choosing a spanning subgraph of $K_{2,2,\dots, 2}$ with independence number $2$, call it $G_1$, and adding a vertex $c$ such that the set of vertices $S =  V(G_1) - N_G[c]$ is an independent $\mathcal{M}$-alternating dominating set in $\overline{G_1}$.

Let $V_i = \{x_i, y_i\}$ be the vertices of each of the color classes of $K_{2,2,\dots, 2}$.  Note that since the set $S$ forms a clique in $G$, $c$ is adjacent to at least one vertex from $V_i$ for each $i \in [k]$.  For convenience, if $c$ is not adjacent to a vertex in a class $V_i$, we will label that vertex $y_i$ (and the other vertex in the set $x_i$). To show that $G$ is well-hued with sequence $(2,4,\dots, 2k, 2k+1)$, we need to consider maximal $\ell$-colorable subgraphs $T$.  To do so, we will induct on $\ell$. We first show that $G$ is well-covered. 
Note that $\alpha(G_1)=2$ so if $T$ is an independent set with $n(T)\ge 3$, then $T$ must contain $c$. However, this contradicts the assumption that $S$ is a clique in $G$. Therefore, there are no independent sets of size $3$ or more. On the other hand, we know that each vertex of $G$ has degree at most $n(G)-2$ so $T$ cannot have order $1$. It follows that $G$ is well-covered. 

Now suppose for some $1 \le \ell' < k+1$ that we know that all maximal $j$-colorable subgraphs have order $2j$ when $j \le \ell'$. Let $\ell=\ell'+1$ and let $T$ be a maximal $\ell$-colorable subgraph. We know that $n(T) < 2\ell+1$ for otherwise $\alpha(G)>2$. So we shall assume that $n(T)<2\ell$. Note that by Theorem~\ref{thm:222}, we know that $G_1$ is well-hued. Therefore, $T$ must contain the vertex $c$. 

Choose an  $\ell$-coloring with color classes $\mathcal{U} = \{U_1, U_2, \dots,  U_{\ell}\}$. Since $n(T)<2\ell$, then we may assume that  $|U_{\ell}| = 1$. Let $T' = U_1\cup U_2\cup\dots \cup U_{\ell-1}$. We can see that $T'$ is an $(\ell-1)$-colorable subgraph of $G$.  If $n(T') < 2(\ell - 1)$ vertices,  then by the induction hypothesis $T'$ is not maximal, and another vertex could be added (which cannot be the vertex in $U_{\ell}$, or $T$ is $(\ell - 1)$-colorable).  Therefore, we may assume that $n(T) = 2\ell - 1$.

Let $\mathcal{M}_{\mathcal{U}}$ be the matching in the complement of $\overline{G}$ corresponding to the classes of $\mathcal{U}$ with two vertices, and let $\mathcal{M}_{\mathcal{V}}$ be the matching in the complement corresponding to $\mathcal{V} = \{V_1, V_2, \dots, V_k\}$.    

Now consider the symmetric difference of $\mathcal{M}_{\mathcal{U}}$ and $\mathcal{M}_{\mathcal{V}}$.  The symmetric difference consists of alternating ($\mathcal{M}_{\mathcal{U}}, \mathcal{M}_{\mathcal{V}}$) paths or even cycles. In fact, by re-coloring the vertices of $T$, we may assume that there are no alternating cycles. We can consider two cases, based on the vertex from $U_{\ell}$. 

We claim that we may assume $c$ is the lone vertex in $U_{\ell}$. Indeed, if the vertex in $U_{\ell}$ is from $G_1$, call it $w_i \in V_i$, then we note that it must be the endpoint of an ($\mathcal{M}_{\mathcal{U}}, \mathcal{M}_{\mathcal{V}}$)-alternating path $P$ where the vertex of $P$ adjacent to $w_i$ is also in $V_i$ as there is no edge in $\mathcal{M}_{\mathcal{U}}$ that saturates $w_i$. Moreover, the other terminal vertex of $P$ must be $c$, in which case we can recolor the vertices along the path so that both vertices in the same partite set in $G_1$ receive the same color and $c$ is the lone vertex in its color class. So we may assume that $c$ is the lone vertex in $U_{\ell}$. Note that we may further assume that every vertex of $S$ is in $V(T)$ for otherwise there exists some $s \in S - V(T)$ such that  $V(T) \cup \{s\}$ is $\ell$-colorable where $s$ is given the color $\ell$. 

Pick maximal $(M_{\mathcal{U}}, M_{\mathcal{V}})$-alternating path $P$ of length at least $2$ and write 
\[P= x_{\alpha_1}y_{\alpha_2}x_{\alpha_2}y_{\alpha_3}\dots x_{\alpha_k}y_{\alpha_{k+1}}\]
such that $V_{\alpha_i} = \{x_{\alpha_i}, y_{\alpha_i}\}$ for $2 \le i \le k$ and $x_{\alpha_i}$ and $y_{\alpha_{i+1}}$ are assigned the same color in $\mathcal{U}$ for $1 \le i \le k$.
Thus, $y_{\alpha_1} \not\in V(T)$ for otherwise $P$ is not maximal. Assume first that for some $y_i \in S$ that $y_i = y_{\alpha_j}$ is on $P$. In this case, if we reindex the color classes in $\mathcal{U}$ so that $U_i = \{x_{\alpha_i}, y_{\alpha_{i+1}}\}$ for $1 \le i \le j-1$, then 
\[\mathcal{U'} = \bigcup_{1\le i \le j-1}V_{\alpha_i} \cup \left(\mathcal{U} - \{U_1, \dots, U_{j-1}, U_{\ell}\}\right) \cup \{c, y_{\alpha_j}\}\]
induces an $\ell$-colorable subgraph $H$ with $V(H) = V(T) \cup \{y_{\alpha_1}\}$. Therefore, we shall assume that no vertex of $S$ is on $P$. 

Since $S$ is a $M_{\mathcal{V}}$-alternating dominating set, there exists a shortest $(M_{\mathcal{V}}, E(\overline{G}) - M_{\mathcal{V}})$-alternating path $Q$ in $\overline{G}$ starting at some $y_{\beta_1} \in S$ and ending at some $y_{\beta_j}$ on $P$. Therefore,  we can write $Q = y_{\beta_1}x_{\beta_1}y_{\beta_2}x_{\beta_2}\dots y_{\alpha_j}$ such that  $y_{\beta_1}\in S$, $y_{\beta_j}$ is on $P$, and $V_{\beta_1} \subset \mathcal{U}$.  Next, note that $y_{\beta_2}$ is in some set in $\mathcal{U}$ for otherwise we could let 
\[\mathcal{U'} = \left(\mathcal{U} - (V_{\beta_1}\cup U_{\ell})\right) \cup \{y_{\beta_1}, c\} \cup \{x_{\beta_1}, y_{\beta_2}\}\]
which induces an $\ell$-colorable graph $H$ with $V(H) = V(T) \cup \{y_{\beta_2}\}$. This in turn implies that $x_{\beta_2}$ is in the same set as $y_{\beta_2}$ in $\mathcal{U}$ for otherwise $y_{\beta_2}$ is the last vertex on $Q$. Continuing this argument, if $y_{\beta_3}$ is not in $T$, then we could define 
\[\mathcal{U'} = \left(\mathcal{U} - (V_{\beta_1}\cup V_{\beta_2} \cup U_{\ell})\right) \cup \{y_{\beta_1}, c\} \cup \{x_{\beta_1}, y_{\beta_2}\}\cup \{x_{\beta_2}, y_{\beta_3}\}\]
which induces an $\ell$-colorable graph $H$ with $V(H) = V(T) \cup \{y_{\beta_3}\}$. Thus, we finally arrive at the conclusion that all vertices of $Q$ are saturated by $\mathcal{U}$ in $\overline{G}$. In this case, since we have assumed $ y_{\alpha_j}$ is on $P = x_{\alpha_1}y_{\alpha_2}x_{\alpha_2}y_{\alpha_3}\dots x_{\alpha_k}y_{\alpha_{k+1}}$ and $y_{\alpha_1}\not\in V(T)$ and all vertices on $Q$ and $P$ are saturated by $\mathcal{U}$ in $\overline{G}$, then the path 
\[R = cy_{\beta_1}x_{\beta_1}y_{\beta_2}x_{\beta_2}\dots y_{\beta_i}x_{\beta_i}y_{\alpha_j}x_{\alpha_{j-1}} y_{\alpha_{j-1}}\cdots x_{\alpha_1}y_{\alpha_1}\]
is such that every vertex other than  $y_{\alpha_1}$ are saturated by $\mathcal{U}$. In this case, we can recolor all vertices of $R$ so that each $\{y_{\alpha_k}, x_{\alpha_k}\}$ receive the same color for $1 \le k \le j-1$, each $\{y_{\beta_k}, x_{\beta_k}\}$ receive the same color for $1 \le k \le i$, and $\{c, y_{\alpha_j}\}$ receive the same color, contradicting the maximality of $T$.

\end{proof}

\begin{corollary}
Every well-hued graph $G_1$ with sequence $(2,4, 6, \dots, 2k)$ occurs as a subgraph of some well-hued graph $G$ with sequence $(2,4, 6, \dots, 2k, 2k+1)$.
\end{corollary}

Next, we turn our attention to well-hued graphs with sequence $(2, 4, \dots, 2k, 2k+1, \dots, n(G))$ where $n(G)\ge 3k$. 

\begin{theorem}\label{thm:3k}  $G$ is a connected well-hued graph with sequence $(2, 4, 6, \dots, 2k, 2k+1, \dots, n(G))$ where $n(G)\ge 3k$  and $k \ge 2$ if and only if $G$ is obtained from a clique $G_1 = K_{n(G)-k}$, and a graph $G_2$ of order $k$ where the following are true:
\begin{enumerate}
\item[1.] for any $u \in V(G)$, the graph induced by $V(G) - N_G[u]$ is a clique, and 
\item[2.] for any set of $k$ vertices in $G_1$, say $X_k$, there exists a perfect matching between $X_k$ and $V(G_2)$ in $\overline{G}$.
\end{enumerate}
\end{theorem}

\begin{proof}
Suppose first that $G$ is well-hued with sequence $(2, 4, 6, \dots, 2k, 2k+1, \dots, n(G))$ where $n(G)\ge 3k$. Let $\ell = n(G) -k$, which is the chromatic number of $G$, and let $c$ be any $\ell$-coloring of $G$ with color classes $V_1, \dots, V_{\ell}$ such that $|V_i| = 2$ for $i \in [k]$ and $|V_i| = 1$ for $k+1 \le i \le \ell$. Write $V_i = \{x_i, y_i\}$ for $i \in [k]$ and $V_i = \{x_i\}$ for $k+1 \le i \le \ell$. Note that $V_{k+1} \cup \cdots \cup V_{\ell}$ induces a clique in $G$ for otherwise $\chi(G) < \ell$. Moreover, if there exists some $i \in [k]$ and pair of vertices $\{x_j, x_t\} \subset V_{k+1} \cup \cdots \cup V_{\ell}$ such that $x_j$ is not adjacent to $x_i$ and $x_t$ is not adjacent to $y_i$, then $V_1 \cup \cdots \cup V_k \cup \{x_j, x_t\}$ is $(k+1)$-colorable with $x_i$ and $x_j$ given the same color and $x_t$ and $y_i$ given the same color. However, this contradicts the assumption that all maximal $(k+1)$-colorable subgraphs of $G$ have order $2k+1$. Therefore, we may assume that if $x_j$ and $x_t$ are not adjacent to some vertex in $V_i$ for $i \in [k]$, then it is the same vertex in $V_i$. That is, relabeling if necessary, we may assume each vertex in $\{x_{k+1}, \dots,  x_{\ell}\}$ is adjacent to every vertex in $\{x_1, \dots, x_k\}$. Now we know there exists a perfect matching in $\overline{G}$ between $\{x_{k+1}, \dots, x_{2k}\}$ and $\{y_1, \dots, y_k\}$ since $G$ is well-hued. Relabeling $x_{k+1}, \dots, x_{2k}$ if necessary, assume that $y_i$ not adjacent to $x_{k+i}$ for $i \in [k]$. If for some $1 \le \alpha < \beta \le k$ $x_{\alpha}$ is not adjacent to $x_{\beta}$, then we could recolor the vertices of $G$ so that $x_{k+\alpha}$ and $y_{\alpha}$ receive the same color, $x_{k+\beta}$ and $y_{\beta}$ receive the same color and $x_{\alpha}$ and $x_{\beta}$ receive the same color, implying $\chi(G)<\ell$ which is a contradiction. Therefore, $\{x_1, \dots, x_k\}$ induces a clique and we have $\{x_1, \dots, x_{\ell}\}$ indeed induces a clique. Moreover, since $G$ is well-hued, Properties (1) and (2) hold. 

In the other direction, suppose that $G$ is obtained from $G_1 = K_{n(G) - k}$ and a graph $G_2$ of order $k$ that satisfies Properties (1) and (2). We show that $G$ is well-hued with sequence $(2, 4, 6, \dots, 2k, 2k+1, \dots, n(G))$. It is clear that $\chi(G) = \omega(G) = n(G) -k$. Note that by Property (1), there is no maximal independent set of size greater than $3$. On the other hand, by Condition (2), each vertex in $G_1$ is not adjacent to some vertex in $G_2$ and vice versa. Therefore, $G$ is well-covered with $\alpha(G) =2$. All that remains is to show for any $2 \le \ell \le n(G) - k -1$, every maximal $\ell$-colorable subgraph has order $2\ell$ if $\ell\le k$ or $\ell + k$ if $\ell >k$.

Let $H$ be any maximal $\ell$-colorable subgraph. If $\ell \le k$, we know that $H$ cannot have order strictly greater than $2\ell$ for this would imply that $H$ contains an independent set of size $3$. On the other hand, if $\ell > k$, then $H$ cannot have order at least $\ell + k +1$ as this would imply that two vertices in $G_1$ have the same color. So assume that $H$ has order at most $2\ell-1$ if $\ell \le k$ and order at most $\ell + k -1$ if $\ell > k$. Partition the vertices of $H$ into color classes  $W_1, \dots, W_{\ell}$ where for some $j \in [\ell]$, $|W_i| = 2$ for $i \in [j]$, and $|W_i| =1$  for $j+1 \le i \le \ell$.  It follows that $W_{j+1} \cup \cdots \cup W_{\ell}$ induces a clique for otherwise $\chi(H)<\ell$. Let $F = V(G) - V(H)$. We may assume that each vertex in $W_{j+1}\cup \cdots \cup W_{\ell}$ is adjacent to each vertex in $F$ for otherwise $H$ is not maximal. Write $W_i = \{w_i\}$ for $j+1 \le i \le \ell$ and $W_i = \{x_i, y_i\}$ for $i \in [j]$ (from above we may assume $j<k$).

Reindexing if necessary, we assume color classes  $W_1, \dots, W_{\alpha}$ contain one vertex from $G_1$ and one vertex from $G_2$ and for $W_{\alpha+1}, \dots, W_j$, each $W_i$ contains both vertices from $G_2$. Further, choose a coloring $c$ of $H$ which maximizes $\alpha$. Label the vertices of $V(G_1) \cap W_i$ as $x_i$ for $i \in [\alpha]$ and $V(G_2) \cap W_i$ as $y_i$ for $i \in[\alpha]$.
\vskip2mm
\noindent\textbf{Case 1:} Suppose $V(G_1) \cap \{w_{j+1}, \dots, w_{\ell}\} \ne \emptyset$. Relabeling if necessary, suppose $w_{j+1} \in V(G_1)$. If $w_{j+1}$ is not adjacent to some $z \in W_{\alpha+1} \cup \cdots \cup W_j$, then we could recolor $z$ with $j+1$ and contradict our choice of $c$. So $w_{j+1}$ is adjacent to every vertex in $W_{\alpha+1} \cup \cdots \cup W_j$. Also, $w_{j+1}$ is adjacent to every vertex in $V(G_2) - V(H)$ for otherwise $H$ is not maximal. By Property (2), $w_{j+1}$ is not adjacent to some vertex $y_i \in \{y_1, \dots, y_{\alpha}\}$. Consider the coloring $c'$ that agrees with $c$ for all vertices except $x_i$ is given color $j+1$ and $w_{j+1}$ is given color $i$. Applying all the same arguments to $x_i$, the lone vertex with color $j+1$, we find that $x_i$ is adjacent to every vertex of $V(G_2) - \{y_1, \dots, y_{i-1}, y_{i+1}, \dots, y_{\alpha}\}$. By Property (2) applied to any $k$-set containing $\{x_i, w_{j+1}\}$, there exists $y_m \in \{y_1, \dots, y_{i-1}, y_{i+1}, \dots, y_{\alpha}\}$ that is either not adjacent to  $x_i$ or not adjacent to $w_{j+1}$. If $w_{j+1}$ is not adjacent to $y_m$, then we could give $w_{j+1}$ the color $m$ and $x_m$ the color $j+1$ and again argue that $x_m$ is adjacent to every vertex in $V(G_2) - \{y_1, \dots, y_{\alpha}\}$. On the other hand, if $x_i$ is not adjacent to $y_m$, then we could assign $w_{j+1}$ the color $i$, $x_i$ the color $m$, and $x_m$ the color $j+1$ and reach the same conclusion. Continuing this same argument we reach the conclusion that either $N_{\overline{G}}(\{x_1, \dots, x_{\alpha}, w_{j+1}\})\subseteq \{y_1, \dots, y_{\alpha}\}$, contradicting Property (2), or $H$ is not maximal. Therefore, this case cannot occur. 

\vskip5mm
\noindent\textbf{Case 2:} Suppose $\{w_{j+1}, \dots, w_{\ell}\}  \subseteq V(G_2)$. We may assume that $w_{j+1}$ is adjacent to each vertex of $V(G_1) - \{x_1, \dots, x_{\alpha}\}$ for otherwise $H$ is not maximal. Therefore, there exists $x_i \in \{x_1, \dots, x_{\alpha}\}$ not adjacent to $w_{j+1}$ by Property (2). Consider the coloring $c'$ which agrees with $c$ for all vertices except $x_i$ and $w_{j+1}$ are given color $j+1$. Applying the same arguments to $y_i$, we see that $y_i$ is adjacent to every vertex in $V(G_1) - \{x_1, \dots, x_{\alpha}\}$. Applying Property (2) to any $k$-set containing $x_i$, there exists $x_m \in \{x_1, \dots, x_{i-1}, x_{i+1}, \dots, x_{\alpha}\}$ that is not adjacent to either $y_i$ or $w_{j+1}$. Continuing this same argument we reach the conclusion that either  $N_{\overline{G}}(\{y_1, \dots, y_{\alpha}, w_{j+1}\})\subseteq \{x_1, \dots, x_{\alpha}\}$, contradicting Property (2), or $H$ is not maximal. This case cannot occur either.
\vskip5mm
Having exhausted all cases, we may conclude no such subgraph $H$ exists and $G$ is in fact well-hued. 

\end{proof}

In summary, Theorems~\ref{thm:2k+1} and ~\ref{thm:3k} give us a characterization of those well-hued graphs with sequence $(2, 4, \dots, 2k, 2k+1, \dots, n(G))$ when $n(G) =2k+1$ or $n(G) \ge 3k$. Therefore, we pose the following problem to complete the characterization of all well-hued graphs with $a_1 =2$.

\begin{problem} Characterize all connected well-hued graphs with sequence $(2, 4, \dots, 2k, 2k+1, \dots, n(G))$ where $2k+1 < n(G) < 3k$. 
\end{problem}

\section{Well-hued Cographs}\label{sec:cographs}

In this section, we use the same terminology and definitions as that found in \cite{GKM-2022, GKM-2020}. Notably, a $k$-chromatic graph $G$ of order $n$ is \emph{well-equi-hued} if it is well-hued and where $a_i = \frac{i}{k}n$ for $i \in \{1, \dots, k\}$. A graph $G$ is said to be \emph{well-k-partite} if every maximal $k$-partite subgraph has the same cardinality, denoted $\alpha_k(G)$. The class of complement reducible graphs, or cographs, is the class of graphs that contains $K_1$ and is closed under the graph operations of disjoint unions and joins.  Precisely, cographs can be defined recursively by the following three rules:
\begin{enumerate}
\item The graph comprised of a single vertex is a cograph.
\item If $G_1$, $G_2$, $\dots$, $G_k$ are cographs, then so is $G_1 \cup G_2 \cup \dots \cup G_k$ (disjoint unions).
\item If $G_1$, $G_2$, $\dots$, $G_k$ are cographs, then so is $G_1 \vee G_2 \vee \dots \vee G_k$ (joins).
\end{enumerate}
Since cographs are defined in terms of joins and unions, first we consider the joins and unions of graphs.

\subsection{Unions and joins of cographs}
To consider when cographs are well-hued, we first note that the disjoint union of well-hued cographs is well-hued. The join of two graphs was previously considered in \cite{GKM-2022}, in which the following theorem was shown:

\begin{theorem}[\cite{GKM-2022}] \label{thm:join}
Consider the join $G = G_1\vee G_2$.  The graph $G$ is well hued if and only if $G_1$ and $G_2$ are well-equi-hued and $\alpha_1(G_1) = \alpha_1(G_2)$.
\end{theorem}

However, it is not difficult to show that we can say slightly more when $G$ is well-hued and arose from the join of two graphs.

\begin{corollary} \label{cor:join}
If the join $G = G_1 \vee G_2$ is a well-hued graph, then $G$ is also well-equi-hued.
\end{corollary}

\begin{proof}
By the previous theorem, we know $G_1$ and $G_2$ must be well-equi-hued. Furthermore,  $\alpha_1(G_1) = \alpha_2(G_2)$ which implies that
\[
\frac{1}{\chi(G_1)}n(G_1) = \frac{1}{\chi(G_2)}n(G_2). 
\]
Furthermore, by the ``componendo" rule for proportions, we note that for positive $a,b,c,d$, if $a/b = c/d$, then $a/b = c/d = (a + c)/(b + d)$.
Therefore, 
\begin{equation} \label{eq:ratios}
\frac{1}{\chi(G_1)}n(G_1) = \frac{1}{\chi(G_2)}n(G_2) =
\frac{1}{\chi(G_1) + \chi(G_2) }\left( n(G_2) + n(G_1)\right).
\end{equation}

  If we take any $k \leq \chi(G) = \chi(G_1) + \chi(G_2)$, we can consider the size of maximal $k$-partite subgraphs. By assumption $G$ is well-hued, so we just need to consider the order of one maximal $k$-partite subgraph. If we write $k = k_1 + k_2$, with $k_1 \leq \chi(G_1)$ and $k_2 \leq \chi(G_2)$, we can choose a maximal $k_1$-partite subgraph from $G_1$ and a maximal $k_2$-partite subgraph from $G_2$.  Then we get an equation for $a_k$ that we can simplify by Equation \ref{eq:ratios}:
\[
a_k = \frac{k_1}{\chi(G_1)}n(G_1) +  \frac{k_2}{\chi(G_2)}n(G_2)
  = \frac{k_1 + k_2}{\chi(G_1) + \chi(G_2) }\left( n(G_2) + n(G_1)\right),
\]
implying that $G$ is well-equi-hued.
\end{proof}

To determine the family of well-hued cographs, it is sufficient to determine those cographs that are well-equi-hued.  Since every connected cograph must be the join of two graphs, then every well-hued cograph is either well-equi-hued, or a disjoint collection of well-equi-hued graphs.  For that reason, we revisit the disjoint union of graphs to determine when the union of two graphs is well-equi-hued. However, first we make the following simple observation:

\begin{observation} \label{obs:well-equi-hued}
For any well-hued graph $G$ with sequence $(a_1, a_2, \dots)$ and any $k$, $1 \leq k \leq \chi(G)$, 
\[a_k \geq \frac{k}{\chi(G)}n(G).\]
\end{observation}

\begin{theorem} \label{thm:union}
Consider the disjoint union $G  = G_1 \cup G_2$.  The graph $G$ is well-equi-hued if and only if $G_1$ and $G_2$ are well-equi-hued and $\chi(G_1) = \chi(G_2)$.
\end{theorem}

\begin{proof}
In one direction, assume that $G$ is well-equi-hued. First, suppose (for contradiction) that at least one of $G_1$ or $G_2$ is not well-hued.  Without loss of generality, suppose that $G_1$ has maximal $k$-colorable subgraphs $S$ and $S'$ such that $n(S) \neq n(S')$.  Then choose any maximal $k$-colorable subgraph $T$ in $G_2$.  The subgraphs $S \cup T$ and $S' \cup T$ are both maximal $k$-colorable subgraphs of $G$, but with different orders, contradicting the fact that $G$ is well-hued. Therefore, let $(a_1, a_2, \dots)$ represent the well-hued sequence for $G$, and $(b_1, b_2, \dots)$ and $(c_1, c_2, \dots)$ the well-hued sequences for $G_1$ and $G_2$, respectively. Note that the same reasoning as above gives  $a_k = b_k + c_k$ for all $k$.

Since $G$ is well-equi-hued, the differences between successive terms in the well-hued sequence are constant, i.e., $a_2 - a_1 = a_3 - a_2 = \dots = a_{\chi(G)} - a_{\chi(G) - 1}$, where  $a_{\ell} - a_{\ell- 1} = (b_{\ell} - b_{\ell- 1}) +  (c_{\ell} - c_{\ell- 1})$ for all $\ell \leq \chi(G) = \max \{\chi(G_1), \chi(G_2)\}$. Since differences between successive terms in any well-hued sequence are non-increasing, this means that the differences for the sequences $b_i$ and $c_i$ must also have constant differences, out to $\chi(G)$, which also implies that $\chi(G_1) = \chi(G_2) = \chi(G)$.  

Lastly,  the fact that $G$ is well-hued and Observation \ref{obs:well-equi-hued} imply the following for $1 \leq k \leq \chi(G)$:
\[
\frac{k}{\chi(G)}n(G_1)  + \frac{k}{\chi(G)}n(G_2) = 
\frac{k}{\chi(G)}n(G) = a_k
= b_k + c_k
\geq \frac{k}{\chi(G)}n(G_1)  + \frac{k}{\chi(G)}n(G_2)
\]
Therefore, it must be the case that $b_k = (k/\chi(G))n(G_1)$ and $c_k = (k/\chi(G))n(G_2)$; thus $G_1$ and $G_2$ are well-equi-hued with same chromatic number.

Next, suppose that $G_1$ and $G_2$ are well-equi-hued with $\chi(G_1) = \chi(G_2)$.
Then $\chi(G) =\chi(G_1) = \chi(G_2)$. Choosing $k$ for $1 \leq k \leq \chi(G)$, we note that any maximal $k$-colorable subgraph of $G$ must arise by choosing maximal $k$-colorable subgraphs from both $G_1$ and $G_2$. Since $G_1$ and $G_2$ are well-hued, this implies that every maximal $k$-colorable subgraph of $G$ must have order $(k/\chi(G))n(G_1) +(k/\chi(G))n(G_2)  = (k/\chi(G))n(G)$, implying that $G$ is well-equi-hued.
\end{proof}

Corneil, Lerchs, and Burlingham \cite{CLB-1981} proved many equivalent characterizations for cographs, including characterizing cographs as the class of $P_4$-free graphs and the comparability graph of a multitree (implying that cographs are perfect graphs). The fact that a cograph $G$ is a perfect graph means that we get $\chi(G) = \omega(G)$, and we can rewrite the previous theorem as follows:

\begin{corollary}  \label{cor:union}
Consider the disjoint union $G  = G_1 \cup G_2$ where $G_1$, $G_2$ are cographs.  The graph $G$ is well-equi-hued if and only if $G_1$ and $G_2$ are well-equi-hued and $\omega(G_1) = \omega(G_2)$.
\end{corollary}

The similarity between Corollary \ref{cor:join} and Corollary \ref{cor:union} allows us to show that the set of well-equi-ued cographs is closed under complementation.

\begin{theorem} \label{thm:cographcomp}
Let $G$ be a cograph.  Then $G$ is well-equi-hued if and only if $\overline{G}$ is well-equi-hued.
\end{theorem}

\begin{proof}
We can proceed by induction on $n(G)$. For $n(G)= 1$, the isolate trivially satisfies this condition. Our induction hypothesis is that for any cograph $G$ such that $n(G)<k$, if $G$ is well-equi-hued, then $\overline{G}$ is well-equi-hued.

If $G$ is formed by a disjoint union of cographs $G_1$ and $G_2$, then $G_1$ and $G_2$ must be well-equi-hued with  $\omega(G_1) = \omega(G_2)$.  In the complement, this implies that $\overline{G_1}$ and $\overline{G_2}$ are well-equi-hued (by the induction hypothesis) and $\alpha_1(G_1) = \alpha_1(G_2)$.  Therefore, the join of $\overline{G_1}$ and $\overline{G_2}$ is well-equi-hued as well.

Similarly, if $G$ is formed by a join of cographs $G_1$ and $G_2$, then $G_1$ and $G_2$ must be well-equi-hued with  $\alpha_1(G_1) = \alpha_1(G_2)$.  In the complement, this implies that $\overline{G_1}$ and $\overline{G_2}$ are well-equi-hued (by the induction hypothesis) and $\omega(G_1) = \omega(G_2)$.  Therefore, the union of $\overline{G_1}$ and $\overline{G_2}$ is well-equi-hued as well.
\end{proof}

\subsection{Cotrees of cographs}

Another idea introduced in \cite{CLB-1981} by Corneil, Lerchs, and Burlingham is the parse tree for a cograph, called the cotree of a cograph.  The cotree represents the way to parse the cograph into its unions and joins, in the following way:
\begin{itemize}
\item The nodes of the cotree that are leaves are labeled with the vertices of its associated cograph.
\item The internal nodes of the tree (i.e., all non-leaf nodes) are labeled with either 0, corresponding to a graph union, or 1, corresponding to a graph join, and alternate such that no two adjacent vertices have the same label.
\item Every internal node has two or more children, unless the associated cograph is a single vertex.
\end{itemize}

\begin{figure}[h!]
\begin{center}
\includegraphics[scale = .4]{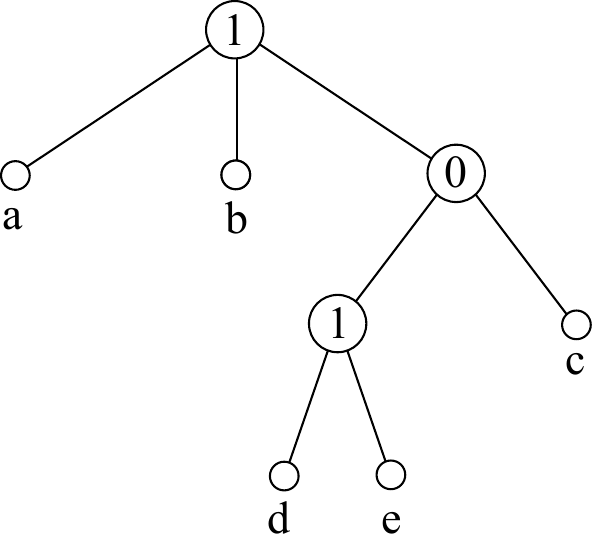} \hspace{50pt}
\includegraphics[scale = .4]{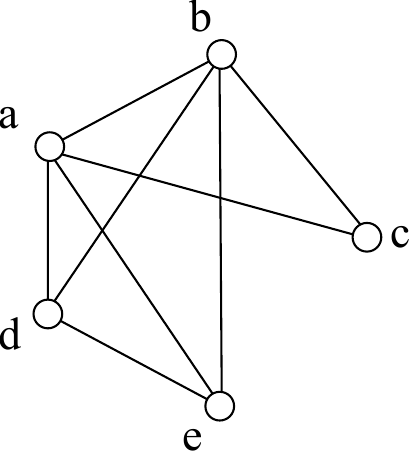}
\caption{A cotree on the left, and the corresponding cograph on the right.\label{fig:cotreeIntro}}
\end{center}
\end{figure}

An example of a cograph, and its associated cotree are given in Figure \ref{fig:cotreeIntro}. There are a few nice properties of cotrees, including the following:
\begin{enumerate}
\item With the description of cotrees given above, cographs can be uniquely mapped to their corresponding cotree, and vice versa.

\item Rooted subtrees of cotrees are cotrees. This follows from the fact that every induced subgraph of a cograph is a cograph. Both of these properties are shown in \cite{CLB-1981}.

\item The cotree of $\overline{G}$ can be found by interchanging the labels on all internal nodes. 
\end{enumerate}

Because of the first property above, we can now focus on the equivalent problem describing the cotrees that correspond to well-equi-hued cographs. In light of Theorem \ref{thm:cographcomp}, all that matters is the underlying tree itself (the labeling of internal nodes doesn't matter). The second property above allows us to use induction, and we can prove our first structural result about such trees.

\begin{theorem} \label{thm:homogeneousChildren}
In the cotree of a well-hued cograph with two or more vertices, the children of an internal node are either all internal nodes, or all leaves.
\end{theorem}

\begin{proof}
We can use induction on the maximum distance from a leaf to the root. As a base case, if the maximum distance is one, then the root is adjacent to only leaves and the result is trivially true.  

Furthermore, we can show that this is the only case for which the root can be adjacent to a leaf. Suppose (for contradiction) that this is not the case. Again, in light of Theorem \ref{thm:cographcomp}, we may assume without loss of generailty that the root is labeled with ``1". We need to show that all other children of the root are leaves.  If this is not the case, there is another child of the root with the label ``0", which contains at least two children. This subtree corresponds to a cograph with independence number at least two. However, the independence number of the cograph corresponding to this subtree, and the single leaf, do not match, which is a contradiction to Theorem \ref{thm:join}.

Now we may assume that any cotree corresponding to 
a well-hued cograph for which the maximum distance from a leaf to the root is $k$ has the property that all children of an internal node are either all internal nodes, or all leaves. Then consider a cotree corresponding to a well-hued cograph for which the maximum distance from a leaf to the root is $k +1$.

Since no neighbor of the root can be a leaf, each child of the root defines a subtree corresponding to a well-equi-hued cograph with two or more vertices, and we can apply the induction hypothesis.
\end{proof}

\begin{figure}
\includegraphics[scale=.3]{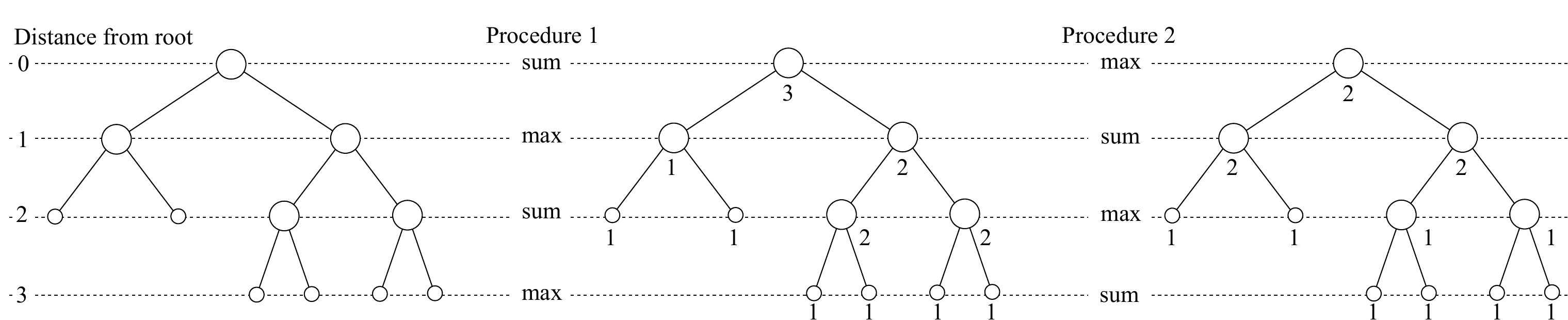}
\caption{An example of a tree satisfying the Uniform Assignment Property \label{fig:cotreeProperty}}

\end{figure}

In addition to recognition of cographs, cotrees are useful algorithmically in that they provide simple algorithms for calculating many graph parameters quickly.  In particular, consider Procedures 1 and 2, illustrated in Figure~\ref{fig:cotreeProperty}.

\begin{table}[!h]
\begin{tabular}{l l}
\begin{minipage}[t]{.45\textwidth}

Procedure 1:
\begin{itemize}
\item Place a value of 1 on each of the leaves
\item At each internal node at an even distance from the root, assign a number equal to the sum of all the values on its children.
\item At each internal node at an odd distance from the root, assign a number equal to the maximum value of all its children.
\end{itemize}

\end{minipage}
&
\begin{minipage}[t]{.45\textwidth}

Procedure 2:
\begin{itemize}
\item Place a value of 1 on each of the leaves
\item At each internal node at an odd distance from the root, assign a number equal to the sum of all the values on its children.
\item At each internal node at an even distance from the root, assign a number equal to the maximum value of all its children.
\end{itemize}

\end{minipage}
\end{tabular}
\end{table}

Originally shown in \cite{CLB-1981}, if the root node is a ``1" node, and Procedure 1 is used, then the label at the root is the value of the largest clique. Similarly, if the root node is a ``0" node, and Procedure 2 is used, then the label at the root is the value of the largest clique. 

To get a necessary and sufficient condition for a cograph $G$ to be well-equi-hued, we can use the procedures described on its corresponding cotree $T_G$.  In particular, we will want the following to hold, which for the sake of brevity, we will call the uniform assignment property.

\vskip5mm
\noindent\textbf{Uniform Assignment Property} \textit{We will say that a cotree $T_G$ of a cograph $G$ satisfies the \textbf{uniform assignment property} if $T_G$, assigned values according to Procedure 1 satisfies the property that all children of a vertex $v$ where $v$ is at odd distance from the root have the same value, and the cotree assigned values according to Procedure 2 satisfies the property that all children of a vertex $v$ where $v$ is at even distance from the root have the same value.}
\vskip5mm

An example of a cotree satisfying the uniform assignment property can be seen in Figure \ref{fig:cotreeProperty}.

\begin{theorem}
A cograph is well-equi-hued if and only if its corresponding cotree satisfies the uniform assignment property.
\end{theorem}

\begin{proof}
We can use induction on the maximum distance from a leaf to the root. As in Theorem \ref{thm:homogeneousChildren}, we note that the base case would be a root adjacent to only leaves. This corresponds to a cograph that is either a clique or a collection of isolated vertices, both of which are well-equi-hued, and it is easy to see that this cotree satisfies the uniform assignment property. Therefore, we may assume that for cotrees corresponding to a cograph for which the maximum distance from a leaf to the root is at most $k$, the corresponding cograph $G$ is well-equi-hued if and only if the cotree satisfies the uniform assignment property.  

Now we want to consider a cograph $G$ with a corresponding cotree $T_G$ so that the maximum distance from a leaf to the root is $k+1$.

In one direction, suppose that $G$ is well-equi-hued, and consider the tree $T_G$.  By Theorem \ref{thm:homogeneousChildren}, none of the children of the root are leaves.  Therefore, let $T_{G_1}, T_{G_2}, \dots, T_{G_k}$ refer to the subtrees rooted at the children of the root.  By Corollaries \ref{cor:join} and \ref{cor:union}, we know that each subtree $T_{G_i}$ must be the cotree of a well-equi-hued cograph $G_i$, and we can apply induction.

Note what Procedure 1 (resp., Procedure 2) run on $T_{G_i}$ (for $1 \le i \le k$) results in the same assigned values (on the vertices of $T_{G_i}$) as when  Procedure 2 (resp., Procedure 1) is run on $T_G$, since the addition of the root changes the parity of the levels. Therefore, to check that the uniform assignment property holds, we just need ensure that the assigned values on the children of the root are all the same when Procedure 2 is used on $T_G$. 
 
 By Theorem \ref{thm:cographcomp}, we may make either choice of label for the root. Without loss of generality, suppose that the root has the label ``0". Then applying Procedure 2 on $T_G$,  restricted to a subtree $T_{G_i}$ ($1 \leq i \leq k$) (now rooted at a ``1" node) would result in a value on the root of $T_{G_i}$ which would be the order of the largest clique of $G_i$.  By Corollary \ref{cor:union}, $\omega(G_1) = \omega(G_2) = \dots = \omega(G_k)$.  Therefore, if $G$ is well-equi-hued, the uniform assignment property holds for $T_G$.

In the other direction, suppose that the uniform assignment property holds for $T_G$.   This implies that either none of the children of the root are leaves, or all of them are. Suppose for contradiction that this is not the case, and the root of $T_G$ has an internal node and a leaf for children.  Label the root of $T_G$ as a ``0" node.  Then let $T_{G_1}$ be a subtree rooted at an internal node that is a child of the root, but is not a leaf.  The root of $T_{G_1}$ must be a ``1" node.  If we consider Procedure 2 on $T_G$, the value at the root of $T_{G_1}$ must be the order of the largest clique of $G_1$. Since the root of $T_{G_1}$ is a ``1" node, this must be at least two.  However, the value on the leaf is 1 which contradicts the fact that $T_G$ satisfies the uniform assignment property.

Therefore, let $T_{G_1}, T_{G_2}, \dots, T_{G_k}$ refer to the subtrees rooted at the children of the root. As before, we note that Procedure 1 (resp. Procedure 2) on $T_G$ when restricted to a subtree $T_{G_i}$ ($1 \leq i \leq k$) results in the same assigned values as Procedure 2 (resp. Procedure 1). Since $T_G$ satisfies the uniform assignment property, so do $T_{G_1}, T_{G_2}, \dots, T_{G_k}$; as a result, $G_1$, $G_2$, $\dots$, $G_k$ are well-equi-hued cographs.

Without loss of gnerality, suppose that the root has the label ``0". Then all children of the root have the same value, and the assigned value at the root of $T_{G_i}$ is the order of the largest clique of $G_i$.  By Corollary \ref{cor:union}, if $\omega(G_1) = \omega(G_2) = \dots = \omega(G_k)$ and $G_i$ are well-equi-hued, then union is well-equi-hued. Therefore, if the uniform assignment property holds for $T_G$, then $G$ is well-equi-hued.
\end{proof}

\section{Concluding Remarks}
In summary, we proved that any connected well-hued graph with $a_2 = a_1 +2$ and $a_1 \ge 3$ is the corona of a complete graph. The next natural question would to characterize all connected well-hued graphs with $a_2 = a_1 + 3$ and $a_1 \ge 4$ (note our example obtained by taking the corona of $K_n$ and attaching a $K_2$ to a leaf). Additionally, does there exist a connected well-hued graph with $a_2 = a_1 + j$ for any $4 \le j \le a_1 - 2$? 

Next, we were able to characterize the connected well-hued graphs with $\alpha(G)=2$ and sequence $(2, 4, 6, \dots, 2k, 2k+1, \dots, n(G))$ provided $n(G) = 2k+1$ or $n(G) \ge 3k$. We leave the remaining cases as an open problem.
\vskip2mm
\noindent\textbf{Problem 1.} \textit{Characterize all connected well-hued graphs with sequence $(2, 4, \dots, 2k, 2k+1, \dots, n(G))$ where $2k+1 < n(G) < 3k$. }
\vskip2mm

Finally, note that our exploration of well-hued cographs stemmed from the graphs in Figure~\ref{fig:table} whose complement was also well-hued. Is there a structural characterization of when both $G$ and $\overline{G}$ are well-hued?


\begin{thebibliography}{99}
\bibitem{CLB-1981}
   D. G.~Corneil, H.~Lerchs, and L. Stewart Burlingham, 
   Complement reducible graphs, 
   Discrete Applied Math., \textbf{3}(3): (1981) 163--174.

\bibitem{FH-1983}
   A. S. Finbow and B. L. Hartnell, 
   A game related to covering by stars,
   Ars Combin., 16: (1983) 189--198.

\bibitem{fhn-1993}
    A. S.~Finbow, B.~L.~Hartnell, and R.~J.Nowakowski,
    A characterization of well-covered graphs of girth {$5$} or greater,
    J. Combin. Theory Ser. B 57 (1993) 44--68.

\bibitem{fhn-1994}
    A. S.~Finbow, B.~L.~Hartnell, and R.~J.~Nowakowski,
    A characterization of well-covered graphs that contain neither {$4$}- nor {$5$}-cycles,
    J. Graph Theory 18 (1994) 713--721.
    
 \bibitem{fhn-1988}
    A. S. Finbow, B. L. Hartnell, and R. J. Nowakoski,
    Well-dominated graphs: a collection of well-covered ones,
    Ars Combin., 25-A: (1988) 5--10.
    



\bibitem{GKM-2022}
   W.~Goddard, K.~Kuenzel, and E.~Melville,
   Well-hued graphs,
   Discrete Applied Math., 320: (2022) 370--380.
   
\bibitem{GKM-2020}
   W.~Goddard, K.~Kuenzel, and E.~Melville,
   Graphs in which all maximal bipartite subgraphs have the same order,
   Aequationes Math., 94: (2020) 1241--1255. 

\bibitem{H-1999}
   B. L. Hartnell, 
   Well-covered graphs,
   J. Combin. Math. Combin. Comput., 29: (1999) 107--115.
   
\bibitem{LPP-1984}
   M. Lesk, M. D. Plummer, and W. R. Pulleyblank,
   Equi-matchable graphs,
   Graph theory and Combinatorics, Academic Press, London, 239-254 (1984)
   
\bibitem{P-1970}
   M.~D.~Plummer, 
   Some covering concepts in graphs, 
   J. Combin. Theory 8: (1970) 91 -- 98.


\bibitem{p-1993}
    M.~D.~Plummer,
    Well-covered graphs: a survey,
    Quaestiones Math. 16 (1993) 253--287.




\end{thebibliography}
\end{document}